\documentclass[runningheads]{llncs}
\usepackage{graphicx}
\usepackage{amsmath}
\usepackage{amssymb}
\usepackage[scr=boondoxo]{mathalpha}
\usepackage{geometry}
\usepackage{tikz}
\usepackage{datetime}
\usetikzlibrary{positioning, arrows.meta, shapes.geometric}
\usepackage{booktabs}

\newcommand\blfootnote[1]{
    \begingroup
    \renewcommand\thefootnote{}\footnote{#1}
    \addtocounter{footnote}{-1}
    \endgroup
}

\begin{document}
\title{An Efficient Integer Programming Model for Solving the Master Planning Problem of
Container Vessel Stowage\thanks{This work is partially sponsored by the Danish Maritime Fund under grant 2021-069.}}
%
%
\titlerunning{An Efficient IP Model for the Master Planning Problem}
\authorrunning{{J. van Twiller et al.}}

\institute{{IT University of Copenhagen,} \\ 
{Rued Langgaards Vej 7, 2300 Copenhagen, Denmark} \\ 
\email{{\{jaiv,rmj\}@itu.dk}}
\and 
{Roskilde University,} \\ 
{Universitetsvej 1, Roskilde, 4000, Denmark} \\ 
\and 
{Aarhus University,} \\ 
{Ny Munkegade 118, 8000 Aarhus, Denmark} \\ 
}

\author{Jaike van Twiller\inst{1} \and
Agnieszka Sivertsen\inst{2} \and 
Rune M. Jensen\inst{1} \and
Kent H. Andersen \inst{3} }
\authorrunning{J. van Twiller et al.}
%
\maketitle              
\begin{abstract}
\vspace{-.75cm}
\blfootnote{\textit{Preprint accepted for 15th International Conference on Computational Logistics}}

A crucial role of container shipping is maximizing container uptake onto vessels, optimizing the efficiency of a fundamental part of the global supply chain. In practice, liner shipping companies include block stowage patterns that ensure that containers in above and below deck partitions of bays have the same destination. Despite preventing restows, increasing free space, and benefits for crane makespan and hydrostatics, this practical planning requirement is rarely included in stowage optimization models. In our paper, we introduce a novel 0-1 IP model that searches in the space of valid paired block stowage patterns, named template planning, which ensures sufficient vessel capacity and limits to crane makespan, trim, and bending moment. Our results show that template planning outperforms traditional allocation planning concerning optimality and runtime efficiency while preserving a sufficiently accurate representation of master planning constraints and objectives. 

\keywords{Integer programming  \and Mathematical modelling \and Computational complexity \and Container stowage planning \and Maritime logistics}
\end{abstract}
\section{Introduction}
Containerized shipping is the backbone of world trade. From 1980 to 2023, the volume of international seaborne trade carried by container vessels grew more than 20-fold from 100 million to 2,200 million tons \cite{statista_container_2024}. It is an environmentally friendly mode of transportation that is politically prioritized \cite{european_commission_freight_2007}. From an operational point of view, however, maximizing the volume of cargo transported by a container vessel is challenging. There are several reasons for this. First, in contrast to a land depot, a container vessel is floating on water and must fulfil complex seaworthiness requirements, including draft, stability, stress forces, and lashing force limits. Second, container vessels sail on closed services between ports and are never empty. Crane moves must be distributed along the vessel such that many cranes can work in parallel to minimize the port stay \cite{jensen_container_2018}. Moreover, cranes can only reach containers from the top of stacks, and minimizing the number of containers that must be restowed to reach containers below them is NP-hard \cite{avriel_container_2000,tierney_complexity_2014}. This latter challenge is a focal point of research in the area (e.g., \cite{pacino_fast_2011,zhu_integer_2020}). In practice, however, the problem is mostly avoided by stowing containers with the same port of discharge (POD) on and below deck \cite{jensen_container_2018}. These so-called {\em paired block stowage} patterns also ensure robustness to uncertain cargo in future ports by clearing as much bottom space as possible. For that reason, the patterns are an operational requirement in practice.\footnote{An exception to this is small feeder vessels that are mostly stowed using the same POD in each stack rather than block patterns.} Despite its significance, there is little previous work on paired block stowage \cite{van_twiller_literature_2023}.  

Due to the dependencies between containers loaded in each service port, maximizing the volume of transported cargo entails solving a multi-port stowage planning optimization problem. A scaleable approach is to decompose the problem into a {\em master planning} and {\em slot planning} problem \cite{van_twiller_literature_2023}. The master planning problem assigns containers to load over a sequence of port calls to bay sections of the vessel. Its main constraints are seaworthiness requirements and sufficient crane utilization. The slot planning problem assigns individual containers to slots in each bay section. Its main constraints are stacking and lashing rules. Experimentally, the master planning problem is the hardest to solve and the most important indicator of loadable volume.

In this paper, we focus on the master planning problem. It has been solved efficiently using MIP models with various sets of constraints \cite{pacino_fast_2011,ambrosino_computational_2015,kebedow_including_2018}, but without using paired block stowage patterns. The limited previous work that uses paired block stowage patterns (e.g., \cite{liu_randomized_2011,pacino_crane_2018,larsen_heuristic_2021}) indicates that they are combinatorially hard. To this end, we contribute a new 0-1 IP formulation of the master planning problem that focuses on paired block stowage named {\em template planning}. In contrast to previous models that allocate cargo to bay sections, we only use decision variables to indicate the port of discharge of each pattern and then require sufficient capacity to load all containers. 

Our experimental evaluation uses data from the representative container vessel stowage planning problem suite, i.e., the largest set of benchmark data publically available to date for representative stowage planning problems \cite{sivertsen_representative_2024}. Our results show that the new formulation outperforms traditional formulations concerning optimality gap and efficiency while preserving a sufficiently accurate representation of master planning constraints and objectives. We also contribute the first complexity result on paired block stowage. As mentioned above, previous work indicates that paired block stowage is combinatorially hard. We reduce the set partitioning problem to the template model, thus showing that searching in paired block stowage patterns is NP-hard.  

The remainder of the paper is organized as follows: Section \ref{sec:domain} introduces the domain of master planning, after which Section \ref{sec:related_work} describes the related work on master planning optimization. In Section \ref{sec:model}, the allocation and template planning models are defined, while Section \ref{sec:results} discusses the computational results of solving both models on benchmark data. Finally, Section \ref{sec:conclusion} concludes the main findings of this study.
\section{Master planning} \label{sec:domain}
Our paper examines explicitly the master planning subproblem of stowage optimization; thus, the main emphasis in this section is on this aspect.  We refer the reader to \cite{jensen_container_2018} for an extensive introduction to stowage planning optimization.

A liner shipping company utilizes a fleet of container vessels that operate on fixed schedules along a closed-loop course, similar to a sea bus service that carries goods between ports rather than people. The complete journey covering all ports is called a \textit{voyage}, typically starting in regions with a surplus supply (e.g., Asia) and navigating towards high-demand discharge areas (e.g., Europe).

Cargo typically comes in three sizes, 20ft, 40ft, and 45ft, with standard dimensions of 8ft in width and either 8ft 6inches in height for \textit{dry cargo} or 9ft 6inches for \textit{high cubes}. Container weights range from 4 to 30 tonnes, depending on their load. Each container has a \textit{port of load} (POL) and \textit{port of discharge} (POD). There can also be different types of cargo that require specific handling, e.g., refrigerate containers \textit{reefers}, that need to be stowed near an electric plug, or \textit{IMDG} containers carrying dangerous goods that are often segregated.

\begin{figure}[h!]
\vspace{-0.5cm}
\centering
\includegraphics[scale=0.4]{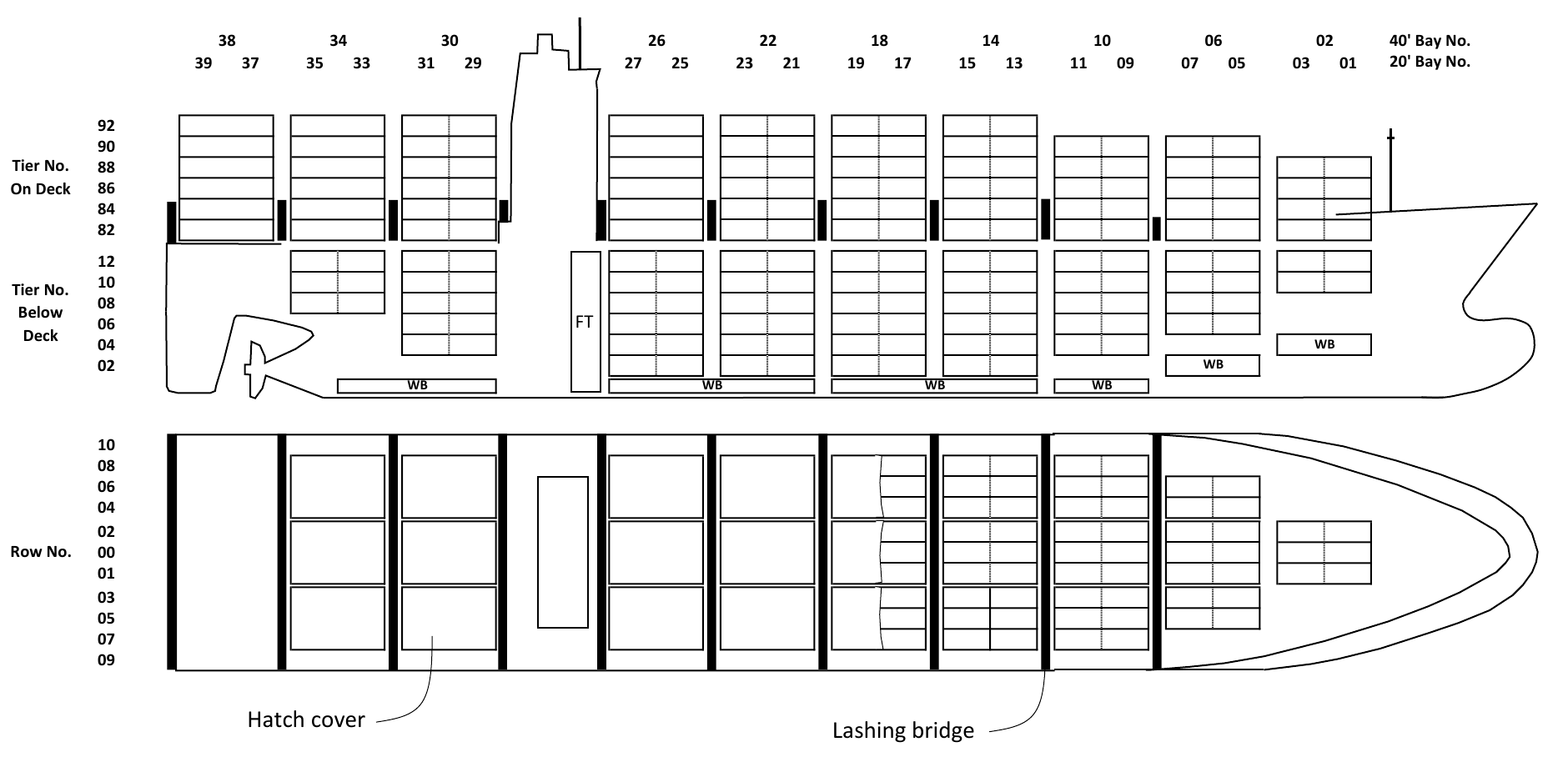}
\caption{Vessel side and top view \cite{van_twiller_literature_2023}.}
\label{fig:vessel1}
\vspace{-0.75cm}
\end{figure}

\begin{figure}[h!]
\vspace{-0.5cm}
\centering
\includegraphics[scale=0.5]{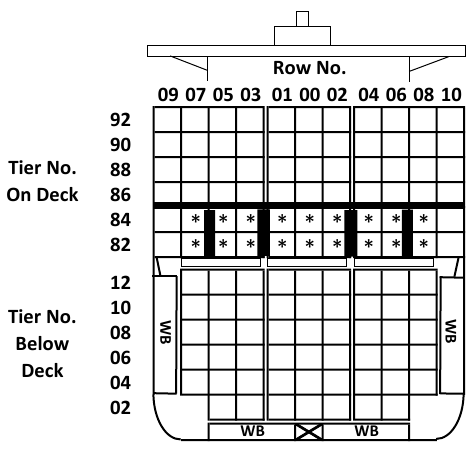}
\caption{Vessel front view \cite{van_twiller_literature_2023}.}
\label{fig:vessel2}
\vspace{-0.5cm}
\end{figure}

The vessel's structure and cellular layout are portrayed in Figure~\ref{fig:vessel1}, demonstrating a longitudinal and top view, and Figure~\ref{fig:vessel2} presents a view of the bay. The vessel is divided into \textit{bays} (02 - 38) that contain \textit{stacks} made up of \textit{cells}. Each cell can hold one 40ft or 45ft container or two 20ft containers in \textit{slots}. To refer to the front and back of the vessel, we use the terms \textit{fore} and \textit{aft}, respectively. Bays are horizontally divided by \textit{hatch covers} into \textit{sectors} either \textit{on deck} or \textit{below deck}. It is worth pointing out that a single bay can have multiple hatch covers, as illustrated in Figure~\ref{fig:vessel1}, where, e.g., bay 30 has three. The sectors can be divided into \textit{partitions} above or below hatch covers; hence, e.g., bay 30 has six partitions. The vessel also has fuel and ballast water tanks, labeled as \textit{FT} and \textit{WB} in Figure~\ref{fig:vessel1}. Along with capacity restrictions, stacks have weight restrictions, and those located below deck also have height boundaries. Figure~\ref{fig:vessel2} shows cells (82-84) with power plugs for refrigerated containers, symbolized by an asterisk.

The primary objective of master planning is to maximize cargo uptake and to minimize operational costs and additional port fees, often associated with reducing the \textit{port stay} of the vessel. Given the size of the vessel and the region it is operating in, a voyage could cover more than ten ports. Moreover, the capacities of modern vessels go beyond accommodating 20,000 Twenty-Foot-Equivalent (TEU) units. The enormity of this task in itself poses a significant challenge.

The master planning problem does not focus on the exact location of a container on a vessel but rather on organizing containers into groups of slots. As such, several container stowage optimization constraints are not considered. The key considerations are capacity and seaworthiness constraints.

Vessels must comply with safety regulations to ensure seaworthiness, which is heavily influenced by the weight distribution on the vessel. The longitudinal (\textit{LCG}), vertical (\textit{VCG}), and transversal (\textit{VCG}) centres of gravity, measured by \textit{trim}, \textit{metacentric height} ({GM}) and \textit{list} respectively, have to be in their specific bounds. Considerable attention must be paid to stress forces, including \textit{shear force}, \textit{bending moment}, and \textit{torsion moment}, to prevent any damage to the vessel.

Minimizing the \textit{makespan} of quay cranes operating simultaneously is crucial for reducing the duration of a port stay. Generally, due to their size, two \textit{quay cranes} cannot operate on two neighboring bays simultaneously. The makespan is determined by the \textit{longest crane}, which refers to the crane with the most extended operational time for loading and unloading containers. Makespan minimization is crucial as it enables maximum utilization of given cranes in a terminal, thereby preventing unnecessary idle time of the equipment.

Moreover, planners utilize specific arrangements known as \textit{stowage patterns} to develop robust plans for cargo to be loaded in future ports. An example is \textit{block stowage}, where stacks divided by the hatch cover consist solely of containers with the same port of discharge (POD). This concept is expanded in \textit{paired block stowage} to guarantee the same POD in \textit{wing blocks} or rows 03-09 and 04-10 as depicted in Figure \ref{fig:vessel2}. On the other hand, \textit{centre blocks} or rows 00-02 can accommodate cargo with the same or another POD than the wing blocks. This approach offers multiple benefits. Firstly, when the entire block is unloaded, it becomes vacant and prepared for a new arrangement of containers. Secondly, the vessel remains balanced due to the even distribution of containers in the wing blocks, preventing tilting. Thirdly,  blocks with a single POD positively impact the long crane, preventing quay cranes from travelling between bays. Lastly, the number of time-consuming hatch cover moves is minimized because all containers are discharged at the same port. Due to these benefits, paired block stowage is an operational requirement for all but the smallest feeder vessels. 

As formerly stated, bays are divided by hatch covers into above and below deck sections and require the removal of the hatch cover to access the cargo stored below deck. The removal process necessitates using quay cranes to rearrange and move the above deck containers that conceal the below deck cargo, in a process known as a \textit{restow} operation. This situation, referred to as \textit{hatch overstowage}, triggers additional, unneeded moves during the vessel's stay in port.
\section{Related work} \label{sec:related_work}
The container vessel stowage planning problem has progressively gained academic attention, emphasizing its significance in the industry \cite{van_twiller_literature_2023}. Existing publications can be categorized into two main types. The first type includes theoretical works that address notable combinatorial problems \cite{avriel_container_2000,tierney_complexity_2014,roberti_decomposition_2018,ding_stowage_2015}. The second type consists of practical studies that focus on heuristic methods. These methods can handle the complex representations of container vessel stowage planning and could be implemented in the industry \cite{pacino_fast_2011,hamedi_containership_2011,liu_randomized_2011}.

One notable approach to the container vessel stowage planning problem is the hierarchical decomposition method  \cite{wilson_container_2001}. This method breaks down the problem into two subproblems: the multi-port master planning problem (MPP) and the slot planning problem (SPP). The former involves assigning groups of containers into blocks on the vessel, while the latter focuses on arranging containers in their slots. The practicality of this approach is evident, as slot planning can be efficiently managed using heuristics, with its impact on runtime being significantly less than that of the master planning problem.

In the following section, we will focus on the characteristics of the MPP problem and the solution methods employed to address it. The most well-known model in the literature is from Pacino et al. \cite{pacino_fast_2011}, which draws inspiration from Ambrosino et al. \cite{ambrosino_stowing_2004}. This model considers the vessel's seaworthiness while aiming to optimize crane work. The decision variables represent the number of containers from a particular group with identical ports of load and discharge to be stowed in a block. Because containers are classified based on their weight, the weight distribution on the vessel can be calculated with a high level of accuracy. Most models found in the literature are inspired by this mathematical formulation \cite{ambrosino_shipping_2018,kebedow_including_2018,bilican_mathematical_2020}. The only alternative method of formulating the master planning problem was suggested by Chao et al. \cite{chao_minimizing_2021}, where a network-flow representation is used to model the problem. However, it does not consider aspects of seaworthiness. The master planning problem has primarily been tackled using mathematical solvers \cite{pacino_fast_2011,pacino_accurate_2012,chao_minimizing_2021}, matheuristics \cite{ambrosino_mip_2015,ambrosino_shipping_2018}, reinforcement learning \cite{van_twiller_towards_2023}, or a mixed approach of exact methods and heuristics \cite{pacino_fast_2013,bilican_mathematical_2020}.

Existing problem formulations focus mostly on minimizing hatch overstowage (e.g., \cite{liu_randomized_2011,pacino_fast_2011,bilican_mathematical_2020}) instead of enforcing paired block stowage constraints. The strategy is to cluster containers heading for the same POD into a single block to avoid unnecessary crane work and create a consolidated free space when the containers are discharged. In the work by Wilson and Roach \cite{wilson_container_2000}, the objective function considers the paired block stowage aspect, aiming to minimize the number of hatches occupied by containers with different ports of discharge (PODs). Liu et al. \cite{liu_randomized_2011} impose a block stowage restriction, allowing only containers with the same POD to be placed in a block. However, these block sizes are smaller than those in conventional block stowage constraints and vary for spaces above and below the hatch cover. In a more theoretical study, Pacino \cite{pacino_crane_2018} demonstrates that including block stowage and crane intensity in container vessel stowage optimization is complex. The provided mathematical formulation could not be solved within a one-hour limit, hence an LNS-based metaheuristic is introduced to find viable solutions.

Relative to the existing body of research, we propose a new 0-1 integer programming formulation for the MPP. This formulation is unique in that it searches within the space of valid paired block stowage patterns, as opposed to the traditional approach of allocating containers to available space. Furthermore, our formulation incorporates a comprehensive set of representative problem features, including vessel capacity for various container types, crane makespan, trim, and bending moment.

\section{Mathematical programming models of the MPP} \label{sec:model}
In previous mathematical models for the MPP, decision variables allocate containers of different types to {partitions} of bays. Recall that a bay with three hatch covers has six partitions, one over and under each hatch cover. This approach does not scale well when we require paired block stowage patterns. Our new approach is to search in the space of valid paired block stowage patterns. If the vessel has three hatch covers, we obtain two storage areas {blocks} per bay: the centre and the wing pair. In the case of four hatch covers, we obtain three blocks: two centres and a wing pair. We refer to this model as template planning.

Recently, there have been several efforts to leverage machine learning in combinatorial optimization \cite{bengio_machine_2021}, especially reinforcement learning seems able to efficiently construct solutions to hard problems (e.g., \cite{hottung_efficient_2022,van_twiller_towards_2023,mirhoseini_graph_2021}). While reinforcement learning has potential, its current limitations in guaranteeing optimality and feasibility for large-scale problems are challenging, as shown in \cite{van_twiller_towards_2023}. Additionally, due to the relative immaturity of stowage planning, we ought to search for novel problem formulations that outperform traditional formulations \cite{van_twiller_literature_2023}. Given the aforementioned considerations, we believe that utilizing mathematical programming in conjunction with well-established solvers is a worthwhile endeavour in this context.

Subsection \ref{sec:definitions} describes sets, parameters, and assumptions to support the mathematical models. Subsequently, we define the MIP model for allocation planning in Subsection \ref{sec:allocation}, while the 0-1 IP model for template planning is defined in Subsection \ref{sec:template}. To show the computational complexity of template planning, we reduce the set partitioning problem to template planning in Subsection \ref{sec:template_np}. 

\subsection{Definitions and assumptions} \label{sec:definitions}
We introduce relevant sets and problem parameters for the MPP in Table \ref{tab:sets} and \ref{tab:params}, respectively. Most definitions speak for themselves, except for the sets and parameters explained in the following subsections. These parameters are extracted from benchmark data, which will be described in Section \ref{sec:results}. Additionally, the following assumptions are made to obtain a simplified version of reality:
\begin{itemize}
\item Cargo only includes 20 ft. and 40 ft. containers, as well as regular and reefer containers. 
\item Each vessel has an arrival condition, represented in the voyage as port 0. Any demand in subsequent ports must be loaded onto the vessel.
\item Loading and discharge times are equal for all ports and types of cargo. 
\item During any voyage, ballast water tanks are constantly half full.
\end{itemize}
\vspace{-0.5cm}

\begin{table}[h!]
\centering
\caption{Sets of the MPP}
\label{tab:sets}
\begin{tabular}{ll}
\toprule
Ports &  $p \in P=\{0,1,2,\ldots\}$ \\
Ports between $i$ and $j$ &  $p \in P^j_i=\{p \in P \mid i \leq p \leq j \}$ \\
Transport pairs &  $(i,j) \in \mathit{TR}=\{(i,j) \in P^2\mid i<j\}$ \\
Onboard transports &  $(i,j) \in \mathit{TR}^{\mathit{OB}}_p=\{(i,j) \in P^2\mid i\leq p, j > p\}$ \\
Discharge transports & $(i,j) \in \mathit{TR}^{\mathit{D}}_p = \{(i,p) \in P^2\mid i < p\}\;$ \\
Load transports & $(i,j) \in \mathit{TR}^{\mathit{L}}_p = \; \{(p,j) \in P^2\mid j > p\}$ \\
Load or discharge transports & $(i,j) \in \mathit{TR}^{M}_p = \mathit{TR}^{\mathit{L}}_p \cup \; \mathit{TR}^{\mathit{D}}_p$ \\
Vessel bays &  $b \in B =  \{1,2,\ldots\}$ \\
Blocks in bay $b$ & $k \in \mathit{BL}_b = \{1,2,\dots, \mathit{HC}_b\mathit{-}1\}$\\
Adjacent bays & $b' \in B' = \{(1,2),(2,3),\ldots,(|B|\mathit{-}1, |B|)\}$ \\
Bays on fore side of bay $b$ & $b' \in B^\mathit{fore}_b = \{1,2,\dots, b\}$ \\
Cargo length &  $l \in \mathit{CL} = \{20',40'\}$ \\
\bottomrule
\vspace{-0.75cm}
\end{tabular}
\end{table}
\vspace{-0.5cm}

\begin{table}[h!]
\centering
\caption{Parameters of the MPP}
\label{tab:params}
\begin{tabular}{ll}
\toprule
Hatch covers per bay (\#) & $\mathit{HC}_b  \; \forall b \in B$ \\
Regular cargo demand (\#) &  ${D}^{l}_{i,j} \; \forall (i,j) \in \mathit{TR}, l \in \mathit{CL}$ \\
Reefer demand (\#) &  $R^l_{i,j} \; \forall (i,j) \in \mathit{TR}, l \in \mathit{CL}$ \\
Volume per container (TEU) &  ${V}_{l} \; \forall l \in \mathit{CL}$ \\
Average container weight for transport $(i,j)$ (tonnes) &  $\Bar{W}_{i,j} \; \forall (i,j) \in \mathit{TR}$ \\
Estimated crane operations (\#) &  $\hat{O}^k \; \forall k \in \mathit{BL}_b, b \in \mathit{B}$ \\
\midrule
Average longitudinal position of bays (meters) &  $L_b \; \forall b \in \mathit{BL}$ \\
Fore longitudinal position of bays (meters) &  $F_b \; \forall b \in \mathit{BL}$ \\
Longitudinal distance between $L_b$ and $F_{b'}$ (meters) &  $\mathit{LD}^{b'}_b \; \forall b \in \mathit{B}, b' \in \mathit{B}$ \\
centre of buoyancy (meters) &  $Z_{p,b} \; \forall p \in \mathit{P}, b \in \mathit{B}$ \\
\midrule
Slot capacity (TEU) &  $C^k_V \; \forall k \in \mathit{BL}_b, b \in \mathit{B}$ \\
Reefer capacity (\#) &  $C^k_R \; \forall k \in \mathit{BL}_b, b \in \mathit{B}$ \\
Weight capacity (tonnes) &  $C^k_W \; \forall k \in \mathit{BL}_b, b \in \mathit{B}$ \\
Maximum crane makespan (\#) &  $\overline{Y}_p \; \forall p \in \mathit{P}$ \\
LCG bounds (meters) &  $\underline{LCG}_p,\overline{LCG}_p \; \forall p \in \mathit{P}$ \\
Bending moment bounds (Newton meters) &  $\underline{\mathit{BM}}_p,\overline{\mathit{BM}}_p \; \forall p \in \mathit{P}$ \\
\bottomrule
\end{tabular}
\vspace{-0.75cm}
\end{table}

\subsection{Allocation planning model} \label{sec:allocation}
Here, we define a typical MIP formulation that allocates cargo to blocks with capacity, crane makespan, hydrostatics, and block stowage constraints inspired by \cite{pacino_crane_2018}. At the core of the allocation model is the minimization of hatch overstowage, which is an NP-hard task \cite{tierney_complexity_2014}. Let $x^k_{i,j} \in \{0,1\}$ indicate whether block $k$ contains cargo of transport $(i,j)$. Let $y^{k,l}_{i,j} \in \mathbb{N}_{0}$ be stowed non-reefer containers and $z^{k,l}_{i,j} \in \mathbb{N}_{0}$ represent stowed reefer containers for transport $(i,j)$, cargo length $l$ and block $k$.

\begin{align}
\text{min } & \displaystyle{\sum_{p \in P}\sum_{(i,j) \in  \mathit{TR}^{\mathit{OB}}_p}  \sum_{b \in B} \sum_{k\in BL_b}  x^k_{i,j}} \label{obj_alloc}
\end{align}
\vspace*{-0.45cm}
{\small\allowdisplaybreaks
\begin{align}
\text{s.t. }        & {\displaystyle \sum_{j \in P^n_{p+1}} x^k_{p,j} \leq 1}
                    & {\displaystyle \forall p \in P^{n-1}_0, k\in \mathit{BL}_b, b \in B} \label{unique_alloc} \\
                    & {\displaystyle \sum_{l \in \mathit{CL}}   y^{k,l}_{i,j} + z^{k,l}_{i,j} \leq M  x^{k}_{i,j}}
                    & {\displaystyle \forall (i,j) \in \mathit{TR}}, k \in \mathit{BL}_b, b \in \mathit{B} \label{link_alloc} \\
                    & {\displaystyle \sum_{b \in \mathit{B}} \sum_{k \in \mathit{BL}_b} y^{k,l}_{i,j} + z^{k,l}_{i,j} = {D}^{l}_{i,j}}
                    & {\displaystyle \forall (i,j) \in \mathit{TR}, l \in \mathit{CL}} \label{volume_alloc} \\
                    & {\displaystyle \sum_{b \in \mathit{B}} \sum_{k \in \mathit{BL}_b} z^{k,l}_{i,j} = R^l_{i,j}}
                    & {\displaystyle \forall (i,j) \in\mathit{TR}, l \in \mathit{CL}} \label{rf_volume_alloc} \\
                    & {\displaystyle \sum_{(i,j) \in \mathit{TR}^{\mathit{OB}}_p} \sum_{l \in \mathit{CL}} {V}_{l}(y^{k,l}_{i,j} +  z^{k,l}_{i,j}) \leq C_V^k \sum_{(i,j) \in \mathit{TR}^{\mathit{OB}}_p} x^k_{i,j}} 
                    & {\displaystyle \forall p \in P, k\in \mathit{BL}_b, b \in \mathit{B}}\label{capacity_alloc} \\
                    & {\displaystyle \sum_{(i,j) \in \mathit{TR}^{\mathit{OB}}_p}  \sum_{l \in \mathit{CL}} z^{k,l}_{i,j} \leq C_R^k \sum_{(i,j) \in \mathit{TR}^{\mathit{OB}}_p} x^k_{i,j}} 
                    & {\displaystyle \forall p \in P, k\in \mathit{BL}_b, b \in \mathit{B}} \label{rf_capacity_alloc} \\
                    & {\displaystyle \sum_{(i,j) \in \mathit{TR}^{\mathit{OB}}_p} \sum_{l \in \mathit{CL}}  \Bar{W}_{i,j}(y^{k,l}_{i,j} +  z^{k,l}_{i,j}) \leq C_W^k \sum_{(i,j) \in \mathit{TR}^{\mathit{OB}}_p} x^k_{i,j}} 
                    & {\displaystyle \forall p \in P, k\in \mathit{BL}_b, b \in \mathit{B}} \label{weight_capacity_alloc} \\
                    & {\displaystyle \sum_{b \in b'}\sum_{k \in \mathit{BL}_{b}} \sum_{(i,j) \in \mathit{TR}^M_p} \sum_{l \in \mathit{CL}} y^{k,l}_{i,j} + z^{k,l}_{i,j} \leq \overline{Y}_p} & {\displaystyle \forall p \in P^{n-1}_1, b' \in B'} \label{makespan_alloc}   \\
                    & {\displaystyle \underline{\mathit{LCG}}_p} \leq {\displaystyle \sum_{b \in B} L_b \sum_{k \in \mathit{BL}_b} \sum_{l \in \mathit{CL}} \sum_{(i,j) \in  \mathit{TR}^{\mathit{OB}}_p} \Bar{W}_{i,j}(y^{k,l}_{i,j} + z^{k,l}_{i,j})} \leq {\displaystyle \overline{\mathit{LCG}}_p} & \forall p \in \mathit{LP}^{n-1}_1 \label{trim_alloc}\\
                    & {\displaystyle \underline{\mathit{BM}}_b} \leq {\displaystyle \sum_{b' \in B^{\mathit{fore}}_b} \mathit{LD}^{b'}_b \sum_{k \in \mathit{BL}_{b'}} \sum_{l \in \mathit{CL}} \sum_{(i,j) \in  \mathit{TR}^{\mathit{OB}}_p} (\Bar{W}_{i,j}(y^{k,l}_{i,j} + z^{k,l}_{i,j})} & \nonumber \\  
                    &  \quad\quad\quad\quad\quad\quad\quad\quad\quad\quad\quad\quad\quad\quad\quad\quad\quad\quad 
                    - Z_{p,b'}) \leq {\displaystyle \overline{\mathit{BM}}_b} & {\displaystyle \forall b \in B, p \in \mathit{LP}^{n-1}_1 \label{bending_alloc}}
\end{align}}

The objective \ref{obj_alloc} expresses the desire to minimize the use of blocks. Constraint \ref{unique_alloc} ensures that blocks have at most one POD during the legs of the voyage. Note that a voyage with $n$ ports has $n-1$ legs. Constraint \ref{link_alloc} links the decision variables, where $x$-variables will equal 1 if the sum of $y,z$-variables are positive. Constraint \ref{volume_alloc} enforces that all cargo in ${D}^{l}_{i,j}$ must be loaded by $y,z$, whereas constraint \ref{rf_volume_alloc} enforces that all reefers in $R^l_{i,j}$ must be loaded by $z$. 

In constraints \ref{capacity_alloc}-\ref{weight_capacity_alloc}, we define block capacity constraints that act as on-off constraints based on onboard $x$-variables. Constraint \ref{capacity_alloc} limits the total TEU per block by the capacity parameter $C^{V}_k$ and parameter ${V}_{l}$ represents the TEU volume of one container with length $l$. Constraint \ref{rf_capacity_alloc} limits reefer utilization by reefer capacity $C^{R}_k$ per block $k$. In this model, a single reefer plug is used for both 20' and 40' reefers because individual cells are not a concern. Constraint \ref{weight_capacity_alloc} limits the total weight of block $k$ by maximum weight $C^{W}_k$ and computes the expected weight of regular and reefers cargo by $\Bar{W}_{i,j}$ as the average weight per TEU during transport $(i,j)$. It is worth noting that the on-off capacities are not strictly necessary, as, e.g., $C^k_V\sum_{(i,j)\in \mathit{TR}^{\mathit{OB}}_p}x_{i,j}^k$ can be replaced by $C^k_V$. Nonetheless, some initial experimenting showed that the on-off capacities shorten the runtime of the allocation model.

Constraint \ref{makespan_alloc} limits the total moves per port with upper bound $\overline{Y}_p$. Let us define the maximum crane makespan as $\overline{Y}_p = \max({Y}_p, \hat{Y}_p)$, where $Y_p$ refers to the maximum long crane provided by data and $\hat{Y}_p=\max_k(C^k_V/1.5)$ is the expected number of container moves of the largest block assuming an equal mix of 20-40 containers. Due to the estimation of crane operations $\hat{O}^k$ in template planning, this formulation is required to ensure feasibility. Details are provided in the next subsection.

In the hydrostatic constraints \ref{trim_alloc} and \ref{bending_alloc}, we quantify over load ports between port $i$ and $j$, which is defined by set $p \in \mathit{LP}^j_i=\{p \in P \mid \sum_{(i,j)\in \mathit{TR}^L_p} \sum_{l \in \mathit{CL}} (D^l_{i,j} + R^l_{i,j}) > 0\}$. Constraint \ref{trim_alloc} sets the upper and lower limit on the longitudinal centre of gravity (LCG) to conform with limits posed on the vessel's trim. These limits, as defined in Equations \ref{for:trim_lb} and \ref{for:trim_ub}, depend on the displacement $d_p$ (the weight of the containers, tanks $W^{\mathit{tw}}$, and lightship $W^\mathit{lsw}$ at departure from port $p$), which could be computed for each instance since it is given that all containers from the load list have to be loaded. Subsequently, the longitudinal centre of buoyancy $\mathit{lcb}(d_p)$ and the trim factor $\mathit{trf}(d_p)$ are interpolated from the hydrostatics table. Additionally, $\underline{t}$ and $\overline{t}$ are the lower and upper bounds of trim, where small instances use $-2.5,2.5$ meters, and the rest use $-2,2$ meters, respectively. $L_b$ is the longitudinal position of bay $b$'s midpoint in meters.
\begin{align}
\underline{LCG}_p = d_p(\mathit{lcb}(d_p) - \frac{\overline{t}}{\mathit{trf}(d_p)}) - \sum_{b \in B}L_b(W^{\mathit{tw}}_b + W^{\mathit{lsw}}_b) \label{for:trim_lb}\\
\overline{LCG}_p = d_p(\mathit{lcb}(d_p) - \frac{\underline{t}}{\mathit{trf}(d_p)}) - \sum_{b \in B}L_b(W^{\mathit{tw}}_b + W^{\mathit{lsw}}_b) \label{for:trim_ub}
\end{align}

Constraint \ref{bending_alloc} sets the upper and lower limit on the bending moment. The bounds are defined by Equations \ref{bm_lb} and \ref{bm_ub}, where $\underline{\mathit{bm}}_b, \overline{\mathit{bm}}_b$ are bounds provided by data, and $W^{\mathit{tw}}, W^\mathit{lsw}$, are the tank and lightship weight. Additionally, $\mathit{LD}^{b'}_b$ is a longitudinal distance between the fore endpoint of bay $b$ and the midpoint of bay $b'$, $Z_{p,b}$ is buoyancy force at departure from port $p$ in bay $b$ interpolated from Bonjean data using displacement $d_p$, and $B^{\mathit{fore}}_b$ is a set of bays positioned on the fore side of bay $b$. The bending moment for bay $b$ is calculated by summing the product of the resulting forces per bay situated on the fore side of $b$ and their respective distances to the fore side point of bay $b$.

\begin{align}
\underline{\mathit{BM}}_b = \underline{\mathit{bm}}_b - {\displaystyle \sum_{b' \in B^{\mathit{fore}}_b} \mathit{LD}^{b'}_b (W^{\mathit{tw}}_b + W^{\mathit{lsw}}_b) } \label{bm_lb}\\
\overline{\mathit{BM}}_b = \overline{\mathit{bm}}_b - {\displaystyle \sum_{b' \in B^{\mathit{fore}}_b} \mathit{LD}^{b'}_b (W^{\mathit{tw}}_b + W^{\mathit{lsw}}_b) } \label{bm_ub}
\end{align}

\subsection{Template planning model} \label{sec:template}
The primary contribution of this paper is the template planning model, which, to our knowledge, has not been previously considered. The main idea is to eliminate all decision variables except the block indicators $x^k_{i,j}$ and ensure that the blocks designated for storing containers provide sufficient capacity. This approach potentially enhances scalability. However, it is also less expressive because it does not specify exactly which and how many containers should be stowed in each block. Despite this limitation, reasonable assumptions can be made to address the issue. 
First, since the objective is to minimize the number of used blocks, we can assume that these blocks are fully utilized when in use. Second, although we cannot model the weight of individual containers, it is realistic to assume that each transport is characterized by a specific weight profile. These assumptions are applied to both the allocation and template planning models to enable direct comparison. With these assumptions, we can sufficiently approximate hydrostatics and crane moves for master planning. To maintain brevity, the definitions provided in the Subsection \ref{sec:allocation} are also applicable in this subsection.

\begin{align}
\text{min } & \displaystyle{\sum_{p \in P}\sum_{(i,j) \in  \mathit{TR}^{\mathit{OB}}_p} \sum_{b \in B} \sum_{k\in BL_b} x^k_{i,j} } \label{obj}
\end{align}
\vspace*{-0.45cm}
{\allowdisplaybreaks
\small\begin{align}
\text{s.t. }        & {\displaystyle \sum_{j \in P^n_{p+1}} x^k_{p,j} \leq 1}
                    & {\displaystyle \forall p \in P^{n-1}_0, k\in \mathit{BL}_b, b \in B} \label{unique} \\
                    & {\displaystyle \sum_{l \in \mathit{CL}}{V}_l ({D}^{l}_{i,j} + R^l_{i,j}) \leq M  \sum_{b \in B} \sum_{k \in \mathit{BL}_b} x^{k}_{i,j}}
                    & {\displaystyle \forall (i,j) \in \mathit{TR}} \label{link} \\
                    & {\displaystyle \sum_{(i,j) \in \mathit{TR}^{\mathit{OB}}_p} \sum_{l \in \mathit{CL}}{V}_l ({D}^{l}_{i,j} + R^l_{i,j}) \leq  \sum_{b \in B} \sum_{k \in \mathit{BL}_b} C_V^k \sum_{(i,j) \in \mathit{TR}^{\mathit{OB}}_p} x^k_{i,j}}        & {\displaystyle \forall p \in P} \label{volume} \\
                    & {\displaystyle \sum_{(i,j) \in \mathit{TR}^{\mathit{OB}}_p} \sum_{l \in \mathit{C}} R^l_{i,j} \leq \sum_{b \in B} \sum_{k \in \mathit{BL}_b} C_R^k \sum_{(i,j) \in \mathit{TR}^{\mathit{OB}}_p} x^k_{i,j}}      & {\displaystyle \forall p \in P} \label{reefer} \\
                    & {\displaystyle \sum_{(i,j) \in \mathit{TR}^{\mathit{OB}}_p} \Bar{W}_{i,j} ( \sum_{l \in \mathit{CL}} {D}^{l}_{i,j} + R^l_{i,j})  \leq \sum_{b \in B} \sum_{k \in \mathit{BL}_b} C_W^k \sum_{(i,j) \in \mathit{TR}^{\mathit{OB}}_p} x^k_{i,j}}       & {\displaystyle \forall p \in P} \label{weight} \\ 
                    & {\displaystyle \sum_{b \in b'} \sum_{k \in \mathit{BL}_b} \sum_{(i,j) \in \mathit{TR}_p^M} \hat{O}^k x^k_{i,j} \leq \overline{Y}_p} & {\displaystyle \forall p \in P^{n-1}_1, b' \in B'} \label{makespan}   \\
                    & {\displaystyle \underline{\mathit{LCG}}_p} \leq {\displaystyle \sum_{b \in B} L_b \sum_{k \in \mathit{BL}_b} \sum_{(i,j) \in \mathit{TR}^{\mathit{OB}}_p} \bar{W}_{i,j}C^{V}_k x^k_{i,j}} \leq {\displaystyle \overline{\mathit{LCG}}_p} & {\displaystyle \forall p \in \mathit{LP}^{n-1}_1} \label{trim} \\
                    & {\displaystyle \underline{\mathit{BM}}_b} \leq {\displaystyle \sum_{b' \in B^{\mathit{fore}}_b} \mathit{LD}^{b'}_b (\sum_{k \in \mathit{BL}_{b'}} \sum_{(i,j) \in \mathit{TR}^{\mathit{OB}}_{p} } \bar{W}_{i,j}C^{V}_k x^k_{i,j}} - Z_{p,b'}) \leq {\displaystyle \overline{\mathit{BM}}_b} & {\displaystyle \forall p \in \mathit{LP}^{n-1}_1,  b \in B} \label{bending}            
\end{align}}

The objective \ref{obj} expresses that we want to use as few blocks as possible. Constraint \ref{unique} expresses that a block at most can be assigned to one POD at a time. Constraint \ref{link} links the $x$-variables to the cargo demand, where at least one $x$-variables among the blocks must be equal to 1 if $\sum_{l \in \mathit{CL}}V_l({D}^l_{i,j} + {R}^l_{i,j})$ is positive for some transport $(i,j)$. Constraint \ref{volume} enforces that there must be enough TEU capacity to fit all onboard demand, constraint \ref{reefer} enforces this for onboard reefers, and constraint \ref{weight} ensures this for onboard weight.

In constraint \ref{makespan}, we must assume the number of crane moves required by some block to approximate crane makespan. We could define an expected containers per transport parameter $\hat{O}_{i,j}= \sum_{l \in \mathit{CL}} (D^l_{i,j} + R^l_{i,j} )/  \sum_{b \in B} \sum_{k \in \mathit{BL}_b} x^k_{i,j} \; \forall (i,j) \in \mathit{TR}^\mathit{M}_p$. However, this causes quadratic terms in Equation \ref{makespan}. Instead, we can estimate both parameters by assuming blocks are fully loaded. The minimization of block usage causes highly utilized blocks, which causes most blocks to be loaded fully. Hence, the crane workload is approximated by $\hat{O}^k = C_{\mathit{V}}^k/1.5$, which obtains the expected number of containers per block $k$ by assuming an equal mix of 20 and 40 ft. cargo. This approximation overestimates the moves in blocks as multiple load moves can be associated with a block. Nonetheless, it provides an estimate of how many crane moves are needed in adjacent bays. 

Constraints \ref{trim} and \ref{bending} have matching interpretations as Constraints \ref{trim_alloc} and \ref{bending_alloc} in the allocation model. Nonetheless, similar to the crane makespan, we must approximate the weight in a certain block to compute hydrostatics. Since the total weight and block usage are known, one way could be to find the average weight per block by $\Bar{W}_{i,j}(\sum_{l \in \mathit{CL}} D^l_{i,j} + R^l_{i,j}) / \sum_{b \in B} \sum_{k \in \mathit{BL}_b} x^k_{i,j} \; \forall (i,j) \in \mathit{TR}^\mathit{OB}_p$. As this also leads to a quadratic term, we again assume that blocks are fully loaded, and therefore block weights are approximated by $\bar{W}_{i,j}C^{V}_k$. Yet again, this overestimates the weight in blocks but also provides an approximation of the vessel's hydrostatics. 

\subsection{Template planning is NP-hard} \label{sec:template_np}
Even though paired block stowage ensures a plan without restows, it is still an NP-hard problem. We prove that the template planning problem is NP-hard by reducing the set partitioning problem to a decision version of it.  Recall that the set partitioning problem is the task of deciding whether a given multiset $S$ of positive integers can be partitioned into two subsets $S_1$ and $S_2$ such that the sum of the numbers in $S_1$ equals the sum of the numbers in $S_2$. We translate a set partitioning problem to a template planning problem as follows. Let the vessel consist of $|S|$ blocks. For each element $s \in S$, there is exactly one block $k$ with a volume capacity $C^k_V$ equal to $s$. There are three ports $P = \{1,2,3\}$. In the first port, there are $\sum_{s \in S} s$ containers to load, each with a volume of one TEU. Half of the containers have POD 2, and the other half have POD 3. There are no containers to load in ports 2 and 3. The containers are assumed to have zero weight, and none are reefers. All lower bounds ($\underline{\mathit{LCG}}_p$,$\underline{\mathit{BM}}_b$) are assumed to be minus infinite. In contrast, all upper bounds ($\overline{LCG}_p$, $\overline{\mathit{BM}}_b$,$\overline{\mathit{Y}}_p$) are assumed to be plus infinite. Hence, constraints \ref{reefer}-\ref{bending} have no effect. Since a block in port 1 can only be assigned to one POD that either is 2 or 3, a feasible solution to this template planning problem divides the blocks into two subsets with equal total capacity. Hence, if the template planning problem has a solution, so does the corresponding set partitioning problem and vice versa. Since the reduction can be done in polynomial time, we have shown that the template planning problem is NP-hard.
\section{Results} \label{sec:results}
In this section, we will provide a computational comparison between the allocation and template models. We use the representative container vessel stowage planning problem suite, the largest set of benchmark data available for representative stowage planning problems, including vessels and problem instances \cite{sivertsen_representative_2024}. A vessel summary is provided in Table \ref{tab:vessel_metrics}, whereas the test instances are summarized in Table \ref{tab:instance_metrics}. Additionally, the vessel data contains the hydrostatics table, from which the longitudinal centre of buoyancy $\mathit{lcb}(d_p)$ and trim factor $\mathit{trf}(d_p)$ can be interpolated based on displacement $d_p$. The vessel data also contains the bonjean data used to interpolate the centre of buoyancy $Z_{p,b'}$ based on displacement $d_p$, as well as the vessel weight parameters, vessel distances, maximum capacities, and hydrostatic limits. The instance data provides cargo and reefer demand, average container weight and crane makespan limits. We refer to the benchmark suite \cite{sivertsen_representative_2024} for a detailed explanation. Furthermore, stowage planners need to respond to changed circumstances quickly, therefore an algorithmic runtime of longer than an hour is hard to use in practice \cite{jensen_container_2018}.

The experiments are run on a Linux machine with AMD EPYC 7742 64-Core Processor and 256 GB memory, running on 2.25GHz/3.4 GHz. The mathematical model is implemented in Python 3.9 and solved with CPLEX 22.1. 

\begin{table}[h!]
\vspace{-0.25cm}
\centering
\caption{Vessel metrics with TEU referring to the total capacity in TEU, Reefers refers to the total reefer capacity in plugs, weight is the total weight capacity in tonnes, bays are the number of bays able to hold cargo and Hatch refers to the maximum number of hatches.}
\begin{tabular}{lrrrrr}
    \toprule
    Vessel Size & TEU  & Reefers & Weight  & Bays   & Hatch \\
    \midrule
    Small & 1,040  & 251 & 22,005 & 8     & 1 \\
    Medium & 6,532  & 1,160 & 162,834   & 18    & 3 \\
    Large & 13,482 & 940 & 274,298   & 22    & 4 \\
    Extra Large & 18,854 & 956 & 342,760    & 24    & 4 \\
    \bottomrule
\end{tabular}
\label{tab:vessel_metrics}
\vspace{-0.75cm}
\end{table}

\begin{table}[h!]
\centering
\caption{Instance metrics with \# being the number of instances, ports being the average loading ports per voyage, cargo being the sum cargo demand in TEU on average, reefers being the sum of reefer demand in TEU on average, and AC being the arrival condition as the sum of onboard cargo in TEU on average. The averages are found by computing the arithmetic mean over instances.}
\begin{tabular}{lrrrrr}
    \toprule
    Vessel Size & \#   & Ports & Cargo & Reefers & AC  \\
    \midrule
    Small & 19  & 8.16 &	1,716 & 59 & 654
     \\
    Medium & 21  &  5.57 & 	5,768 &	155 & 4,890 
    \\
    Large & 16 & 3.06 & 7,153 & 211  & 6,785 \\
    Extra Large & 13 & 5.31 & 17,461  & 248 & 11,222  \\
    \bottomrule
\end{tabular}
\label{tab:instance_metrics}
\vspace{-0.25cm}
\end{table}

Table \ref{tab:results} provides a computational comparison between the allocation and template model applied to the instances mentioned above. It should be mentioned that the solver times out at 3,600 seconds, and then returns the objective, gap and runtime of the solution with the best optimality gap found.
The optimality gap of the allocation model is consistently larger than that of the template model across all instance sets, highlighting the difficulty the allocation model faces under these specific constraints. Especially in the larger instances, the template model can find near-optimal solutions while the allocation struggles to find solutions with an optimality gap close to the 0\% optimum. 
Moreover, the template model is significantly faster than the allocation model regarding runtime, with an average speed-up for an instance set ranging from 2 to 4.5 times. Furthermore, a steep increase in computational time is observed from small to medium, large, or extra-large instances, showcasing the complexity of solving real-life instances. Consequently, we argue that the template model scales well to industrial-sized instances, which is not the case for the allocation model.

\begin{table}[h!]
\vspace{-0.25cm}
\centering 
\caption{Comparison of allocation and template models on several sets of instances grouped by vessel size with \# number of instances. Obj. represents the objective value, the duality gap is denoted by Gap (\%), and Time (s) represents the runtime in seconds. The solver either accepts solutions with a duality gap of 1\% or returns the best solution after a runtime limit of 3,600 seconds. All metrics are averaged across instances with the arithmetic mean.}
\begin{tabular}{@{\hspace{0.2cm}}l@{\hspace{0.2cm}}l@{\hspace{0.2cm}}|@{\hspace{0.2cm}}r@{\hspace{0.2cm}}r@{\hspace{0.2cm}}r@{\hspace{0.2cm}}|@{\hspace{0.2cm}}r@{\hspace{0.2cm}}r@{\hspace{0.2cm}}r}
\toprule
& &  \multicolumn{3}{c|@{\hspace{0.2cm}}}{Allocation} & \multicolumn{3}{c}{Template} \\
Vessel Size & \# & Obj. & Gap (\%) & Time (s) & Obj. & Gap (\%) & Time (s)   \\
\midrule
Small & 19 & 88.68 & 0.21 &  5.88 & \textbf{80.80} & \textbf{0.05} & \textbf{1.82} \\
Medium & 21 & 176.90 & 8.49 & 3,371.93 &  \textbf{149.05} & \textbf{0.91} & \textbf{1,387.92} \\
Large & 16 & 169.50 & 9.91 &  3,620.65 & \textbf{163.00} & \textbf{0.99} & \textbf{802.58} \\
Extra Large & 13 & 288.54 & 12.71 &  3,614.61 & \textbf{248.69} & \textbf{0.97} & \textbf{1,654.78} \\
\bottomrule
\end{tabular}
\label{tab:results}
\vspace{-0.25cm}
\end{table}

From Table \ref{tab:results}, one may observe that the objective values between the allocation and template model are different. We can assume a linear relationship between the objective and optimality gap to compute the expected optimal objective $\mathbb{E}[\mathit{Obj}^*]= (1-\mathit{Gap})\mathit{Obj}$. Table \ref{tab:results_obj} shows the difference between the $\mathbb{E}[\mathit{Obj}^*]$ of both models across each set of instances, which shows a mean absolute error $\mu_\mathit{AE}$ of at most 10\% relative to either $\mathbb{E}[\mathit{Obj}^*]$. Additionally, the standard deviation $\sigma_\mathit{AE}$ indicates a reasonable variation in the absolute error, as the coefficient of variation $\mathit{CV}_\mathit{AE}= \sigma_\mathit{AE}/\mu_\mathit{AE}$ remains below 1. Hence, this suggests different but similar expected optimal objective values across instances. The difference, however, is mainly due to the allocation model having more freedom to assign cargo to various blocks and the approximations for the long crane and hydrostatics in template planning, imposing slightly different constraints on the optimization problem. Nonetheless, both models adhere to paired block stowage patterns and aim to minimize block use. Hence, the objective values are within close range.

\begin{table}[h!]
\centering 
\caption{
Comparing the expected optimal objective of allocation and template models on several sets of instances grouped by vessel size with \# number of instances. The expected optimal objective is represented by $\mathbb{E}[\mathit{Obj}^*]= (1-\mathit{Gap})\mathit{Obj}$, absolute error (AE) refers to absolute difference between $\mathbb{E}[\mathit{Obj}^*]$ of both models with arithmetic mean $\mu_{\mathit{AE}}$, standard deviation $\sigma_{\mathit{AE}}$ and coefficient of variation $\mathit{CV}_{\mathit{AE}}$ over instances.}
\vspace{0.cm}
\begin{tabular}{l@{\hspace{0.2cm}}l@{\hspace{0.2cm}}|r@{\hspace{0.2cm}}|r@{\hspace{0.2cm}}|r@{\hspace{0.2cm}}r@{\hspace{0.2cm}}r@{\hspace{0.2cm}}}
\toprule
& &  \multicolumn{1}{@{\hspace{0.2cm}}c@{\hspace{0.2cm}}|}{Allocation} & \multicolumn{1}{@{\hspace{0.2cm}}c@{\hspace{0.2cm}}|}{Template} & \multicolumn{3}{@{\hspace{0.2cm}}c@{\hspace{0.2cm}}}{Absolute error} \\
Vessel Size & \# & $\mathbb{E}[\mathit{Obj}^*]$ & $\mathbb{E}[\mathit{Obj}^*]$ & $\mu_\mathit{AE}$ & $\sigma_{\mathit{AE}}$ & $\mathit{CV}_{\mathit{AE}}$  \\
\midrule
Small & 19 & 88.45 & 80.75 & 7.30 & 5.38 & 0.738 \\
Medium & 21 & 160.92 & 147.63 & 13.29 & 6.31 & 0.475 \\
Large & 16 & 153.41 & 161.31 & 9.89 & 8.05 & 0.814\\
Extra Large & 13 & 249.47 &  246.26 & 3.82 & 3.19 & 0.834 \\
\bottomrule
\end{tabular}
\label{tab:results_obj}
\end{table}

In conclusion, the template model outperforms the allocation model with respect to optimality and runtime at the cost of approximating the long crane and hydrostatics. This trade-off allows us to solve the IP model with paired block stowage patterns, capacity constraints, maximum crane makespan, trim, and bending moment.

\section{Conclusion} \label{sec:conclusion}
This paper introduces a new 0-1 IP model called template planning to solve the master planning problem with paired block stowage patterns and constraints to limit capacity, crane makespan, trim, and bending moment. In previous work, MIP models allocate containers of different types to blocks on the vessel, which does not scale well if block stowage patterns are included. Instead, our so-called template planning model searches in the space of valid paired block stowage patterns.  

The experiments utilize the latest and largest benchmark suite for representative stowage planning problems. Our findings indicate that the template formulation outperforms the allocation model regarding the optimality gap and runtime while preserving an adequate representation of master planning constraints and objectives. Particularly in the larger instance sets, the allocation model struggles to find near-optimal solutions within an hour of runtime, whereas the template model demonstrates scalability by finding near-optimal solutions with an average speed-up between 2 and 4.5 times per set of instances. Additionally, we reduce the set partitioning problem to the template planning model, showing that the computational complexity of searching in valid paired block stowage patterns is NP-hard. 

In future work, we aim to improve the computational efficiency of both models by improving the mathematical formulation and leveraging the cutting-planes approach. Moreover, we will enhance the approximations of crane makespan and hydrostatics in the template model, thus minimizing differences between allocation and template planning.

%
%

\end{document}